\documentclass[11pt]{article}
\usepackage{amsmath}
\usepackage{mathrsfs}
\usepackage{amsfonts}

\usepackage{amsfonts,amsmath,amssymb}
\usepackage{amssymb,amsfonts,amsmath,
latexsym,epsfig,cite,psfrag,eepic,colordvi}
\usepackage{amscd,graphics}
\newcommand{\qed}{\hfill $\Box $}
\newcommand{\pf}{\noindent {\bf Proof.} }

\newtheorem{theorem}{Theorem}[section]

\begin{document}

\title{Note on parity factors of regular graphs
\thanks{This work was supported by the National Natural Science Foundation of China (No. 11101329)}}

\author{
Hongliang Lu\thanks{Corresponding email:
luhongliang215@sina.com (H. Lu)}
\\ {\small  Department of Mathematics}
\\ {\small Xi'an Jiaotong University, Xi'an 710049, PR China}
}

\date{}

\maketitle

\date{}

\maketitle

\begin{abstract}
In this paper, we obtain a sufficient condition for the existence of
parity factors in a regular graph in terms of edge-connectivity.
Moreover, we also show that our condition is sharp.

\end{abstract}

\section{Preliminaries}

Let $G = (V,E)$ be a graph with vertex set $V (G)$ and edge set
$E(G)$. The number of vertices of a graph $G$ is called the
\emph{order} of $G$ and is denoted by $n$. On the other hand, the
number of edges of $G$ is called the \emph{size} of $G$ and is
denoted  by $e$. For a vertex $v$ of  graph $G$, the number of edges
of $G$ incident with $v$ is called the \emph{degree} of $v$ in $G$
and is denoted by $d_{G}(v)$.
%Throughout this paper, $G$ denotes a simple graph of \emph{order}
%$n$ (the number of vertices) and \emph{size} $e$ (the number of
%edges). We denote vertex set and edge set of $G$ by $V(G)$ and
%$E(G)$, respectively. Let $d_{G}(v)$ denote the degree of a vertex
%$v$ in $G$.
For two subsets $S,T\subseteq V(G)$, let $e_{G}(S,T)$ denote the
number of edges of $G$ joining $S$ to $T$.

Let therefore $g, f : V\rightarrow Z^+$ such that $g(v)\leq f (v)$
and $g(v)\equiv f (v) \pmod 2$  for every $v\in V$. Then a spanning
subgraph $F$ of $G$ is called a $(g, f )$-parity- factor, if
$g(v)\leq d_F(v)\leq f (v)$ and $d_F(v)\equiv f (v) \pmod 2$ for all
$v\in V$. Let $a,b$ be two integers   such that $1\leq a\leq b$ and
$a\equiv b\pmod 2$. If $g(v)\equiv a$ and $f(v)\equiv b$ for all
$v\in V(G)$, then a $(g, f )$-parity- factor is called
\emph{$(a,b)$-parity factor}. When $a=1$, $(a,b)$-parity factor is
called \emph{$(1,n)$-odd factor}.

For a general graph $G$ and an integer $k$, a spanning subgraph $F$
such that $$d_{F}(x)=k \ \mbox{ for all } x\in V(G)$$ is called a
\emph{$k$-factor}.  In fact, a $k$-factor is also a $(k,k)$-parity
factor.

Now let us recall one of the most classic results due to Petersen.

\begin{theorem}[Petersen \cite{Petersen}]
Let $r$ and $k$ be integers such that $1\leq k\leq r$. Every
$2r$-regular graph  has a $2k$-factor.
\end{theorem}

By the edge-connectivity, Gallai \cite{Gallai} proved the following
result.
\begin{theorem}[Gallai \cite{Gallai}]\label{Gallai}
Let $r$ and $k$ be integers such that $1\leq k<r$, and $G$  an
$m$-edge-connected $r$-regular  graph, where $m\geq 1$. If one of
the following conditions holds, then $G$ has a $k$-factor.
\begin{itemize}
\item[$($i$)$] $r$ is even, $k$ is odd, $|G|$ is even, and  $\frac{r}{m}\leq k\leq r(1-\frac{1}{m})$;

\item[$($ii$)$] $r$ is odd, $k$ is even and $2\leq k\leq r(1-\frac{1}{m})$;

\item[$($iii$)$]$r$ and $k$ are both odd and $\frac{r}{m}\leq k$.
\end{itemize}

\end{theorem}

Bollob\'as, Satio and Wormald \cite{Bollobas} improved above result.

\begin{theorem}[Bollob\'as, Saito and Wormald \label{Bollobas}
]\label{Bol} Let $r$ and $k$ be integers such that $1\leq k<r$, and
$G$ be an $m$-edge-connected $r$-regular  graph, where $m\geq 1$ is
a positive integer. Let $m^*\in \{m,m+1\}$  such that $m^*\equiv
1\pmod 2$. If one of the the following conditions holds, then $G$
has a $k$-factor.
\begin{itemize}
\item[$($i$)$] $r$ is odd, $k$ is even and $2\leq k\leq r(1-\frac{1}{m^*})$;

\item[$($ii$)$]$r$ and $k$ are both odd and $\frac{r}{m^*}\leq k$.
\end{itemize}
\end{theorem}

In this paper,  we extend Gallai as well as Bollob\'as, Satio and
Wormald result to $(a,b)$-parity factor. The main tool in our proofs
is the famous theorem of Lov\'asz (see\cite{Lovasz}).

\begin{theorem}[Lov\'asz \cite{Lovasz}]\label{Lovasz}
$G$ has a $(g, f )$-parity factor  if and only if for all disjoint
subsets $S$ and $T$ of $V(G)$,
\begin{align*}
\delta(S,T)&=f(S)+\sum_{x\in T}d_{G}(x)-g(T) -e_{G}(S,T)-\tau\geq 0,
\end{align*}
where $\tau$ denotes the number of components $C$, called $f$-odd
components of $G-(S\cup T)$ such that $e_{G}(V(C),T)+f(V(C))\equiv 1
\pmod 2$. Moreover, $\delta(S,T)\equiv f(V(G)) \pmod 2$.
\end{theorem}

\section{Main Theorem}

\begin{theorem}\label{main}\label{main1}
Let $a,b$ and $r$ be integers such that $1\leq a\leq b<r$ and
$a\equiv b\pmod 2$. Let $G$ be a $m$-edge-connected $r$-regular
graph with $n$ vertices. If one of the following conditions holds,
then $G$ has a $(a,b)$-parity factor.

\begin{itemize}
\item[$($i$)$] $r$ is even, $a,b $ are odd, $|G|$ is even,    $\frac{r}{m}\leq b$ and $a\leq r(1-\frac{1}{m})$;

%\item[(ii)] $r$ is odd, $k$ is even and $2\leq k\leq r(1-\frac{1}{m})$;
%
%\item[(iii)] $r$ and $k$ are both odd and $\frac{r}{m}\leq k$.

\item[$($ii$)$]  $r$ is odd, $a,b$ are even and $a\leq
r(1-\frac{1}{m^*})$;

\item[$($iii$)$]  $r$, $a,b$ are odd and $\frac{r}{m^*}\leq b$.
\end{itemize}
\end{theorem}

\pf By Theorem \ref{Bollobas}, (ii) and (iii) are followed. Now we
prove (i). Let $\theta_1=\frac{a}{r}$ and $\theta_2=\frac{b}{r}$.
Then $0<\theta_1\leq \theta_2<1$. Suppose that $G$ contains no
$(a,b)$-parity factors. By Theorem \ref{Lovasz}, there exist two
disjoint subsets  $S$ and $T$ of $V(G)$ such that $S\cup T\neq
\emptyset$, and
\begin{align}\label{eq4}
-2\geq \delta(S,T)&=b|S|+\sum_{x\in T}d_{G}(x)-a|T|
-e_{G}(S,T)-\tau,
\end{align}
where $\tau$ is the number of $a$-odd (i.e. $b$-odd) components  $C$
of $G-(S\cup T)$. Let $C_{1},\cdots,C_{\tau}$
 denote $a$-odd components of $G-S-T$  and $D=C_{1}\cup \cdots\cup C_{\tau}$.

 Note that
\begin{align*}
-2\geq \delta(S,T)&=b|S|+\sum_{x\in T}d_{G}(x)-a|T|
-e_{G}(S,T)-\tau\\
&=b|S|+(r-a)|T|-e_{G}(S,T)-\tau\\
&=\theta_2 r|S|+(1-\theta_1)r|T|-e_{G}(S,T)-\tau\\
&=\theta_2 \sum_{x\in S}d_{G}(x)+(1-\theta_1)\sum_{x\in T}d_{G}(x)-e_{G}(S,T)-\tau\\
&\geq \theta_2 (e_{G}(S,T)+\sum_{i=1}^\tau
e_{G}(S,C_{i}))+(1-\theta_1)(e_{G}(S,T)+\sum_{i=1}^\tau
e_{G}(T,C_{i}))-e_{G}(S,T)-\tau\\
&=\sum_{i=1}^\tau(\theta_2
e_{G}(S,C_{i})+(1-\theta_1)e_{G}(T,C_{i})-1)+(\theta_2-\theta_1)e_G(S,T)\\
&\geq \sum_{i=1}^\tau(\theta_2
e_{G}(S,C_{i})+(1-\theta_1)e_{G}(T,C_{i})-1).
\end{align*}
Since $G$ is connected and $0<\theta_1\leq \theta_2<1$, so $\theta_2
e_{G}(S,C_{i})+(1-\theta_1)e_{G}(T,C_{i})>0$  for each $C_{i}$.
Hence we will obtain a contradiction by showing that  for every
$C=C_{i}$, $1\leq i\leq \tau$, we have
\begin{align}\label{eq1}
\theta_2 e_{G}(S,C)+(1-\theta_1)e_{G}(T,C)\geq 1.
\end{align}
These inequalities together with the previous inequality imply
\begin{align*}
-2\geq \delta_{G}(S,T)&\geq \sum_{i=1}^{\tau}(\theta_2
e_{G}(S,C_{i})+(1-\theta_1)e_{G}(T,C_{i})-1)\\
&> \sum_{i=1}^{\tau-2}(\theta_2
e_{G}(S,C_{i})+(1-\theta_1)e_{G}(T,C_{i})-1)-2\geq -2,
\end{align*}
which is impossible.
 Since $C$ is a $a$-odd component of $G-(S\cup T)$, we have
\begin{align}\label{eq2}
a|C|+e_{G}(T,C)\equiv 1\ \mbox{(mod 2)}.
\end{align}
 Moreover, since
$$r|C|=\sum_{x\in V(C)}d_{G}(x)=e_{G}(S\cup T, C)+2|E(C)|,$$
we have
\begin{align}\label{eq3}
r|C|=e_{G}(S\cup T,C)\ \mbox{(mod 2)}.
\end{align}
 It is obvious that the two
inequalities $e_{G}(S,C)\geq 1$ and $e_{G}(T,C)\geq 1$ imply
$$\theta_2 e_{G}(S,C)+(1-\theta_1)e_{G}(T,C)\geq \theta_2+1-\theta_1=1.$$
Hence we may assume $e_{G}(S,C)=0$ or $e_{G}(T,C)=0$.

Firstly, we consider (i). If
 $e_{G}(S,C)=0$, then $e_{G}(T,C)\geq m$. Since $a\leq r(1-\frac{1}{m})$, then
$\theta_1\leq 1-\frac{1}{m}$ and so $1\leq (1-\theta_1)m$. By
substituting $e_{G}(T,C)\geq m$ and $e_{G}(S,C)=0$ into (\ref{eq1}),
we have
$$(1-\theta_1)e_{G}(T,C)\geq (1-\theta_1)m\geq 1.$$
%If bote $r$ and $k$ are odd, then by (\ref{eq2}) and (\ref{eq3}), we
%have
%$$|C|+e_{G}(T,C)\equiv 1\ \mbox{(mod 2) and }|C|\equiv e_{G}(T,C)\ \mbox{(mod 2)}.$$
%This is again a contradiction.
If  $e_{G}(T,C)=0$, then $e_{G}(S,C)\geq m$.   %By (\ref{eq2}), we
%have $k|C|\equiv 1$ (mod 2), which implies that $k$ must be odd.
%Hence we may assume that $k$ is odd.
Since $\frac{r}{m}\leq b$, hence $\theta_2 m\geq 1$, and so  we
obtain
$$\theta_2 e_{G}(S,C)\geq \theta_2 m\geq 1.$$
Consequently,  condition (i) guarantees (\ref{eq1}) holds and thus
(i) is true. Consequently the proof is complete. \qed

%In order to prove that (ii) implies Claim 2, it suffices to show
%that (\ref{eq1}) holds under the assumption that $e_{G}(S,C)$ or
%$e_{G}(T,C)=0$. %We first consider the first part of (ii), i.e., we
%%assume that $r$ is odd, $k$ is even and $k\leq r(1-\frac{1}{m^*})$.
%If $e_{G}(S,C)=0$, then by (\ref{eq2}), we have $e_{G}(T,C)\equiv 1$
%(mod 2). Hence $e_{G}(T,C)\geq m^*$, and thus
%$$(1-\theta_1)e_{G}(T,C)\geq (1-\theta_1)m^*\geq 1.$$
%If $e_{G}(T,C)=0$, then by (\ref{eq2}), we have $a|C|\equiv 1$ (mod
%2), which contradicts the assumption that $a$ is even.
%
%We next consider  (iii), i.e., we assume that both $r$ and $b$ are
%odd and $\frac{r}{m^*}\leq b$. If $e_{G}(S,C)=0$, then by
%(\ref{eq2}) and (\ref{eq3}), we have
%\begin{align*}
%|C|+e_{G}(T,C)\equiv 1\ \mbox{(mod 2) and }|C|\equiv e_{G}(T,C)\
%\mbox{(mod 2)}.
%\end{align*}
%This is a contradiction. If $e_{G}(T,C)=0$, then by (\ref{eq2}) and
%(\ref{eq3}), we have
%\begin{align*}
%|C|\equiv 1\ \mbox{(mod 2) and }|C|\equiv e_{G}(S,C)\ \mbox{(mod
%2)},
%\end{align*}
%which implies $e_{G}(S,C)\geq m^*$. Thus
%$$\theta_2 e_{G}(S,C)\geq \theta_2 m^*\geq 1.$$
%Consequently the proof is complete. \qed

\noindent\textbf{Remark:} The edge-connectivity conditions  in
Theorem
\ref{main}are sharp.

We give the description for (i). For (ii) and
(iii), the constructions are similar but slightly more complicated.
Let $r\geq 2$ be an even integer, $a,b\geq 1$ two odd integers and $
2\leq m\leq r-2$ an even integer such that $b < r/m$ or
$r(1-\frac{1}{m})<a$. Since $G$ has a $(a,b)$-parity factor if and
only if $G$ has a $(r-b,r-a)$-parity factor, so we can assume $b <
r/m$. Let $J(r, m)$ be the complete graph $K_{r+1}$ from which a
matching of size $m/2$ is deleted. Take $r$ disjoint copies of $J(r,
m)$. Add $m$ new vertices and connect each of these vertices to a
vertex of degree $r-1$ of $J(r, m)$. This gives an
$m$-edge-connected $r$-regular graph denoted by $G$. Let $S$ denote
the set of $m$ new vertices and $T=\emptyset$.  Let $\tau$ denote
the number of components $C$, called $a$-odd components of $G-(S\cup
T)$ such that $e_{G}(V(C),T)+a|C|\equiv 1 \pmod 2$. Then we have
$\tau=r$, and
$$\delta(S,T)=b|S|+\sum_{x\in T}d_{G-S}(x)-a|T| -\tau(S,T)=bm-r<0.$$
So by  Theorem \ref{Lovasz}, $G$ contains no $(a,b)$-parity factors.

%Since $G$ has a $(a,b)$-parity factor if and only if $G$ has a
%$(r-a,r-b)$-parity factor,    the descriptions for (ii) and (iii)
%are similar. Now we the description for (iii). Let $ m\leq r-2$ be
%an integer and $m\in \{m,m+1\}$   an odd integer. Let $1 \leq a\leq
%b< r$ be three odd integers and  such that $b < r/m^*$. We denote
% a cycle of length $(r+2)$ by $C_{r+2}$ and choose a matching of $C_{r+2}$ $M_{(r-m^*)/2}$ be a
%matching with order $(r+2-m^*)/2$  denoted by $M_{(r+2-m^*)/2}$. Set
%$H(r,m^*)=C_{r+2}- M_{(r+2-m^*)/2}$ and
%$I(r,m^*)=\overline{H(r,m^*)}$. Add $m^*$ new vertices and connect
%each of these vertices to a vertex of degree $r-1$ of $I(r, m^*)$.
%This gives an $m^*$-edge-connected $r$-regular graph denoted by $F$.
%Similarly, let $S$ denote the set of $m^*$ new vertices and
%$T=\emptyset$. Let $\tau$ denote the number of $a$-odd components
%$C$ of $G-(S\cup T)$. Then we have $\tau=r$, and
%$$\delta(S,T)=b|S|+\sum_{x\in T}d_{G-S}(x)-a|T| -\tau(S,T)=bm^*-r<0.$$
%So by  Theorem \ref{Lovasz}, $G$ contains no $(a,b)$-parity factors.

%\begin{theorem}
%Let $1\leq a\leq b<r$ be three odd integers. $G$ contains a
%$(a,b)$-parity factor if and only if $G$ contains a $b$-factors.
%
%\end{theorem}
%
%\begin{theorem}
%Let $a,b$ be even and $r$ odd such that $1\leq a<b\leq r$. $G$
%contains a $(a,b)$-parity factor if and only if $G$ contains a
%$a$-factors.
%
%\end{theorem}

\end{document}